\documentclass[12pt]{amsart}

\usepackage{amssymb}
\usepackage{amsmath}
\evensidemargin\oddsidemargin

\begin{document}

\setcounter{page}{1}

\newtheorem{PROP}{Proposition}
\newtheorem{REMS}{Remark}
\newtheorem{LEM}{Lemma}
\newtheorem{THE}{Theorem}
\newtheorem{DEFI}{Definition\!\!}
\newtheorem{THEA}{Theorem A\!\!}
\newtheorem{THEB}{Theorem B\!\!}
\newtheorem{THEC}{Theorem C\!\!}

\renewcommand{\theDEFI}{}
\renewcommand{\theTHEA}{}
\renewcommand{\theTHEB}{}
\renewcommand{\theTHEC}{}

\newcommand{\eqnsection}{
\renewcommand{\theequation}{\thesection.\arabic{equation}}
    \makeatletter
    \csname  @addtoreset\endcsname{equation}{section}
    \makeatother}
\eqnsection

\def\a{\alpha}
\def\bt{\beta}
\def\B{{\bf B}} 
\def\CC{{\mathbb{C}}} 
\def\cia{c_{\a, \infty}}
\def\coa{c_{\a, 0}}
\def\cua{c_{\a, u}}
\def\cL{{\mathcal{L}}} 
\def\Ea{E_\a}
\def\EE{{\mathbb{E}}} 
\def\fa{f_\a}
\def\Ga{\Gamma}
\def\G{{\bf \Ga}} 
\def\i{{\rm i}}
\def\I{{\bf I}}
\def\K{{\bf K}}
\def\Ka{{\bf K}_\a}
\def\L{{\bf L}}
\def\lbd{\lambda}
\def\lcr{\left[}
\def\lpa{\left(}
\def\lva{\left|}
\def\M{{\bf M}}
\def\NN{{\mathbb{N}}}
\def\pb{{\mathbb{P}}}
\def\rl{{\mathbb{R}}}
\def\prst{{\leq_{st}}}
\def\prost{{\prec_{st}}}
\def\prcvx{{\prec_{cx}}}
\def\rpa{\right)}
\def\rcr{\right]}
\def\rva{\right|}
\def\S{{\bf S}}
\def\X{{\bf X}}
\def\Sa{S_\a}
\def\U{{\bf U}}
\def\Un{{\bf 1}}
\def\cN{{\mathcal{N}}} 
\def\Y{{\bf Y}}
\def\Z{{\bf Z}}
\def\Za{{\bf Z}_\a}

\def\claw{\stackrel{d}{\longrightarrow}}
\def\elaw{\stackrel{d}{=}}
\def\qed{\hfill$\square$}
                  
\title[Comparing Fr\'echet and positive stable laws]
      {Comparing Fr\'echet and positive stable laws}

\author[Thomas Simon]{Thomas Simon}

\address{Laboratoire Paul Painlev\'e, Universit\'e Lille 1, Cit\'e Scientifique, F-59655 Villeneuve d'Ascq Cedex. {\em Email} : {\tt simon@math.univ-lille1.fr}}

\address{Laboratoire de physique th\'eorique et mod\`eles statistiques, Universit\'e  Paris Sud, B\^atiment 100, 15 rue Georges Cl\'emenceau, F-91405 Orsay Cedex.}



\begin{abstract} Let $\L$ be the unit exponential random variable and $\Za$ the standard positive $\a$-stable random variable. We prove that $\{(1-\a)\a^{\gamma_\a} {\bf Z}_\a^{-\gamma_\a}, 0< \a <1\}$ is decreasing for the optimal stochastic order and that $\{(1-\a){\bf Z}_\a^{-\gamma_\a}, 0< \a < 1\}$ is increasing for the convex order, with $\gamma_\a = \a/(1-\a).$ We also show that $\{\Gamma(1+\a) \Za^{-\a}, 1/2\le \a \le 1\}$ is decreasing for the convex order, that $\Za^{-\a}\,\prost\, \Gamma(1-\a) \L$ and that $\Gamma(1+\a)\Za^{-\a} \,\prcvx\,\L.$ This allows one to compare $\Za$ with the two extremal Fr\'echet distributions corresponding to the behaviour of its density at zero and at infinity. We also discuss the applications of these bounds to the strange behaviour of the median of $\Za$ and $\Za^{-\a}$ and to some uniform estimates on the classical Mittag-Leffler function. Along the way, we obtain a canonical factorization of $\Za$ for $\a$ rational in terms of Beta random variables. The latter extends to the one-sided branches of real strictly stable densities.

\end{abstract}

\maketitle
 
\section{Introduction and statement of the main results}

This paper deals with two classical random variables. The first one is the 
positive $\a-$stable random variable $\Za$ ($0<\a<1$) which is defined through its Laplace transform
\begin{equation}
\label{Polar1}
\EE\lcr e^{-\lbd \Za}\rcr \; =\; e^{-\lbd^\a}, \qquad \lbd \ge 0.
\end{equation}
Recall that the density of $\Za$ is not explicit except in the case $\a = 1/2$ - see e.g. Example 2.13 in \cite{S} - where it equals
$$\frac{1}{2\sqrt{\pi x^3}}e^{-\frac{1}{4x}}\,\Un_{(0,+\infty)} (x),$$ 
or in the cases $\a=1/3$ resp. $\a=2/3$ where it can be expressed in terms of a Macdonald resp. a Whittaker function - see Formul\ae\, (2.8.31) and (2.8.33) in \cite{Z1}. The second one is the Fr\'echet random variable with shape parameter $\gamma > 0,$ which is defined as the negative power transformation $\L^{-\gamma}$ of the unit exponential random variable $\L.$ Its density is 
$$\frac{1}{\gamma}\, x^{-(1+\frac{1}{\gamma})}\, e^{-x^{-\frac{1}{\gamma}}}\Un_{(0,\infty)}(x)$$
but its Laplace transform is not given in closed form, except for $\gamma =1$ where it can be expressed as a Macdonald function - see e.g. Example 34.13 in \cite{S}. The importance of these two laws stems from limit theorems for positive independent random variables and we refer to Chapter 8.3 resp. 8.13 in \cite{BGT} for more on this topic. Both positive stable and Fr\'echet distributions are known as "power laws", which means that their survival functions have a polynomial decay at infinity. Other common visual features of these laws are the exponential behaviour of the distribution function at zero and the unimodality. Plotting a Fr\'echet density yields a curve whose shape barely differs from that of a positive stable density, as reproduced in \cite{Z1} pp. 144-145.

It is hence a natural question to investigate how these two laws resemble one another. In the present paper, we propose to compare them in a stochastic sense. We will use the general notion of stochastic ordering, referring to \cite{SS} for an account. If $X, Y$ are two real random variables, we say that $Y$ dominates $X$ for the stochastic order and we write $X\,\prst Y$ if 
\begin{equation}
\label{cxx}
\EE[\varphi(X)]\; \le\; \EE[\varphi(Y)]
\end{equation}
for all increasing functions such that the expectations exist. This is equivalent to the fact that $\pb[X\ge x]\le \pb[Y\ge x]$ for all $x\in \rl$ - see Chapter 1 pp. 3-4 in \cite{SS}. When $X\,\prst Y$ it is possible that $X+c\, \prst Y$ for some $c > 0$ and one can ask for an optimal stochastic order between two random variables, that is that no such $c$ exists. In the framework of positive random variables, let us  introduce the following natural definition.

\begin{DEFI} Let $X, Y$ be positive random variables. We say that $Y$ dominates $X$ for the optimal stochastic order and we write 
$$X\;\prost\; Y$$ 
if $X\,\prst Y$ and if there is no $c > 1$ such that $cX\, \prst Y.$
\end{DEFI}

Our second ordering is more related to the dispersion of random variables. If $X, Y$ are two real random variables, we say that $Y$ dominates $X$ for the convex order and we write 
$$X\;\prcvx\; Y$$ 
if (\ref{cxx}) holds for all convex functions such that the expectations exist. When $X$ and $Y$ have finite first moments with $\EE[X] =\EE[Y],$ then $X\prcvx Y$ is equivalent to the fact that 
$$\int_x^{+\infty}\pb[X\ge y]\, dy \;\le\; \int_x^{+\infty} \pb[Y\ge y]\, dy$$ 
for all $x\in \rl$ - see Chapter 2 p.56 in \cite{SS}. When $X$ and $Y$ have finite first moments, the condition $X\prcvx Y$ entails (by choosing $\varphi(x) = x$ and $\varphi(x)=-x$ above) the normalization $\EE[X] =\EE[Y]$, so that contrary to the stochastic order there is no requirement of an optimal formulation for the convex order which is optimal in itself. We refer to the first part of the book \cite{SS} for more details on stochastic and convex orders, and also for other types of orderings.

To state our results, we need some further notation. Using the terminology of \cite{Z1} p.13, set $\S = e^\X$ where $\X=Y(1,-1,-1,1)$ is the spectrally negative $1-$stable random variable with drift coefficient $-1$ and scale parameter $1,$ viz. the random variable with characteristic function
$$\EE[e^{\i \lbd \X}]\; =\; e^{\i \lbd(\log\vert\lbd\vert -1) -\pi\vert\lbd\vert/2}\; =\; \lpa \frac{\i\lbd}{e}\rpa^{\i\lbd}, \quad \lbd\in\rl.$$ 
The random variable $\S$ is an example of a log-stable distribution. Our first main result exhibits two complete orderings for the random variables $\Z_\a,$ after a certain power transformation.

\begin{THEA} For every $0 < \bt<\a <1$ one has
\begin{equation}
\label{stZ}
\S\;\prost\;(1-\a)\a^{\frac{\a}{1-\a}} {\bf Z}_\a^{\frac{-\a}{1-\a}}\;\prost\; (1-\bt)\bt^{\frac{\bt}{1-\bt}}{\bf Z}_\bt^{\frac{-\bt}{1-\bt}}\;\prost\; \L
\end{equation}
and
\begin{equation}
\label{cxZ}
\L\;\prcvx\;(1-\bt) {\bf Z}_\bt^{\frac{-\bt}{1-\bt}}\;\prcvx\; (1-\a){\bf Z}_\a^{\frac{-\a}{1-\a}} \;\prcvx\; e\S.
\end{equation}
\end{THEA}

\medskip

Our next result is expressed in terms of the Mittag-Leffler random variable of order $\a\in (0,1),$ which is defined by
$$\M_\a\; \elaw\; \Za^{-\a}.$$
The denomination comes from the fact that the Laplace transform of $\M_\a$ is expressed in terms of the Mittag-Leffler function - see (\ref{Polar2}) infra. By the Darling-Kac theory, Mittag-Leffler random variables and the associated stochastic processes appear as limit objects for occupation times of Markov processes - see \cite{Bi1} and the whole Section 8.11 in \cite{BGT}. The random variable $\M_\a$ is also distributed as the first-passage time of an $\a-$stable subordinator - see Theorem 2 in \cite{St}. On the other hand, when $\a\in [1/2, 1),$ the random variable $\M_\a$ has the same law as the running supremum of a certain stable process. More precisely, if $\{X_t, \, t\ge 0\}$ stands for the spectrally negative strictly $(1/\a)-$stable L\'evy process normalized such that
$$\EE[e^{\lbd X_t}]\; =\; e^{t\lbd^{1/\a}}, \qquad \lbd\in \rl^+,$$ 
then it is well-known - see e.g. Example 46.7 in \cite{S} and the whole Section 8  in \cite{Bi2} for more on this topic - that
$${\overline X}_1\; =\; \sup\{X_t, \; t\le 1\}\; \elaw\; \M_\a.$$
The following theorem displays a convex ordering for the above suprema, after suitable normalization.
 
\begin{THEB} For every $1/2 \le \bt<\a <1$ one has
\begin{equation}
\label{cxZa}
\Ga(1+\a)\M_\a\;\prcvx\;\Ga(1+\bt) \M_\bt.
\end{equation}
\end{THEB}
According to Theorem 2.A.3 in \cite{SS}, this result entails that on some probability space there exists a martingale $\{M_t, \;t\in [0,1]\}$ such that
$$M_t\; \elaw\; \Ga(2- t/2)\M_{1-\frac{t}{2}}$$
for every $t\in [0,1].$ Notice that this martingale starts at $\M_1 \elaw\Un$ and finishes at $\Ga(3/2)\M_{1/2},$ which is half-Gaussian. Notice that from (\ref{cxZ}) one can also exhibit another martingale starting at $\L$, finishing at $e\S$ and having 
$$(1-\a)\Za^{\frac{-\a}{1-\a}}\; \elaw\;(1-\a)\M_\a^{\frac{1}{1-\a}}$$ 
as marginal law. Our third and last main result compares the Mittag-Leffler distribution and the exponential law for both stochastic and convex orders. It is noticeable that contrary to Theorem A, the orderings are in the same direction.

\begin{THEC} For every $\a \in(0,1)$ one has
\begin{equation}
\label{cxstL}
\M_\a\;\prost\; \Gamma(1-\a)\L\qquad\mbox{and}\qquad \Gamma(1+\a)\M_\a \;\prcvx\;\L.
\end{equation}
\end{THEC}
 
\medskip
 
The stochastic orderings (\ref{stZ}) and (\ref{cxstL}) entail immediately the following optimal comparisons between the random variable $\Za$ and two Fr\'echet distributions, which motivates the title of the present paper:
\begin{equation}
\label{title}
\a (1-\a)^{\frac{1-\a}{\a}}\L^{-\frac{1-\a}{\a}}\;\prost\;\Z_\a\qquad\mbox{and}\qquad\Gamma(1-\a)^{-\frac{1}{\a}}\L^{-\frac{1}{\a}}\;\prost\;\Z_\a.
\end{equation}
It is interesting to notice that these two Fr\'echet distributions are extremal as far as their possible comparison with $\Za$ is concerned. Indeed, the exact behaviours of the distribution function of $\Za$ at zero and infinity, which will be recalled in (\ref{zero}) and (\ref{ZAL}) below, entail that there is no $\gamma\in [1/\a-1, 1/\a]^c$ and no $\kappa > 0$ such that $\kappa\L^{-\gamma}\prst\, \Za.$ These behaviours also show that there is no $\gamma, \kappa > 0$ such that $\Za\,\prst\,\kappa\L^{-\gamma}.$ On the other hand, for any $\gamma\in (1/\a-1, 1/\a)$ one can prove that there exist some $\kappa > 0$ such that $\kappa\L^{-\gamma}\prst\, \Za.$ But it seems difficult to find a formula for the optimal $\kappa$, even in the explicit case $\a =1/2.$ In the latter case, direct computations show indeed that this amounts to finding the infimum of those $c >0$ such that the equation
$$e^{-x^{\frac{2}{\gamma}}}\; =\; \frac{1}{\sqrt{\pi}}\int_{cx}^\infty e^{-\frac{t^2}{4}} dt$$
has no solution on $(0,+\infty).$ This is an ill-posed problem except for $\gamma =1$ (with $c_{\min}=2$) or $\gamma = 2$ (with $c_{\min}=\sqrt{\pi}$), those two cases which were already handled in Theorems A and C. Hence, at the level of the stochastic order, it seems that (\ref{stZ}) and (\ref{cxstL}) is the best that can be said for the comparison between Fr\'echet and positive stable laws. At the level of the convex order, it follows from (\ref{cxZ}) and (\ref{cxstL}) that
$$ \Ga(1+\a)\Z_\a^{-\a}\;\prcvx\;\L\;\prcvx\;(1-\a)\Z_\a^{-\frac{\a}{1-\a}}.$$
This gives some further information on the relationship between $\Za$ and its two associated extremal Fr\'echet distributions.

The paper is organized as follows. In Section 2 we derive two factorizations of $\Za$
with $\a$ rational in terms of Beta and Gamma random variables, which will play some r\^ole in the proof of Theorem C. These factorizations are interesting in themselves and in Section 7 we extend them to the one-sided branches of all real strictly stable densities, with the help of Zolotarev's duality. In Section 3 we derive some explicit computations on Kanter's random variable, which appears to be the key-multiplicative factor of $\Za$ for this kind of questions. The three main theorems are proved in Section 4 and 5, with two concrete applications which are the matter of Section 6. First, we derive some explicit bounds on the median of $\Za,$ showing that the latter behaves quite differently according as $\a\to 0$ or $\a \to 1.$ Second, we prove some uniform estimates on the classical Mittag-Leffler function, answering an open problem recently formulated in \cite{M}. 

\smallskip
 
\noindent
{\bf Convention and notations.} Throughout the paper, the product of two random variables is always meant as an independent product. We also the make the convention that a product over a void set is equal to 1. We will use repeatedly the Legendre-Gauss multiplication formula for the Gamma function:
\begin{equation}
\label{LG}
(2\pi)^{(p-1)/2}p^{1/2-z}\Ga(z)\; =\; \Ga(z/p)\times\cdots\times \Ga((z+p-1)/p)
\end{equation}
for all $z>0$ and $p\in\NN^*.$ 

\section{Factorizing one-sided stable densities}

In this section we aim at factorizing $\Za$ with $\a$ rational in terms of the Beta random variable $\B_{a,b}$ and the Gamma random variable $\G_c,$ whose respective densities are
$$\frac{\Gamma(a+b)}{\Gamma(a)\Gamma(b)} \, x^{a-1} (1-x)^{b-1} \Un_{(0,1)}(x)\qquad\mbox{and}\qquad\frac{x^{c-1} e^{-x}}{\Gamma(c)} \Un_{(0,\infty)}(x).$$ 
Recall the standard formul\ae\, for the fractional moments:
\begin{equation}
\label{SF}
\EE[\B_{a,b}^s]\; =\;\frac{\Ga(a+s)\Ga(a+b)}{\Gamma(a)\Gamma(a+b+s)} \qquad\mbox{and}\qquad \EE[\G_c^s]\; =\; \frac{\Ga(c+s)}{\Gamma(c)}
\end{equation}
over the respective domains of definition. If $n >p \ge 1$ are two integers, let us define the following indices: $q_0 = 0, q_p = n$ and if $p\ge 2,$  
$$q_j \; =\; \sup\{i\ge 1, \; ip < jn\}$$ 
for all $j = 1, \ldots p-1.$ Notice that the family $\{q_j, \, 0\le j\le p\}$ is increasing in $[0,n].$ Observe also that if $q_{j+1}\ge q_j +2,$ then for all $i\in[q_j+1, q_{j+1}-1]$ one has
$$\frac{i-j}{n-p} - \frac{i}{n}\; =\; \frac{pi-jn}{n(n-p)}\; > \; 0.$$
Last, it is easy to see that
\begin{equation}
\label{O1}
\{i-j, \; j\in [0,p], \, i\in[q_j+1, q_{j+1}-1]\} \; =\; \{1, \ldots, n-p\}
\end{equation}
and that the set on the left-hand side is strictly increasing with respect to the  lexicographic order in $(j,i).$ 

\begin{THE} With the above notation, for all $n > p\ge 1$ integers one has the identities in law
\begin{equation}
\label{BetaGamma}
{\bf Z}_{\frac{p}{n}}^{-p} \;\elaw \; \frac{n^n}{p^p}\;\, \prod_{j=0}^{p-1} \lpa \prod_{i=q_j+1}^{q_{j+1}-1} \G_{\frac{i}{n}}\rpa\; \times\; \prod_{j=1}^{p-1}\; \B_{\frac{q_j}{n},\frac{j}{p} - \frac{q_j}{n}}
\end{equation}
and
\begin{equation}
\label{Beta}
{\bf Z}_{\frac{p}{n}}^{-p} \;\elaw \; \frac{n^n}{p^p (n-p)^{n-p}}\; \L^{n-p}\; \times\; \prod_{j=0}^{p-1} \lpa \prod_{i=q_j+1}^{q_{j+1}-1} \B_{\frac{i}{n},\frac{i-j}{n-p} - \frac{i}{n}}\rpa\; \times\; \prod_{j=1}^{p-1}\; \B_{\frac{q_j}{n},\frac{j}{p} - \frac{q_j}{n}}.
\end{equation}
\end{THE} 

\proof We first show that the two random variables on the right-hand sides of (\ref{BetaGamma}) and (\ref{Beta}) have the same law. Formula (\ref{LG}) and a fractional-moment identification entail
\begin{equation}
\label{IDL}
\L^{n-p}\; \elaw \; (n-p)^{n-p}\; \prod_{k=1}^{n-p} \G_{\!\!\frac{k}{n-p}}\;\elaw \; (n-p)^{n-p} \;\prod_{j=0}^{p-1} \lpa \prod_{i=q_j+1}^{q_{j+1}-1} \G_{\!\frac{i-j}{n-p}}\rpa,
\end{equation}
where the second identity follows from (\ref{O1}). Hence, the random variable on the right-hand sides of (\ref{Beta}) can be written
$$\frac{n^n}{p^p}\,\prod_{j=0}^{p-1} \lpa \prod_{i=q_j+1}^{q_{j+1}-1} \B_{\frac{i}{n},\frac{i-j}{n-p} - \frac{i}{n}}\,\times\,\G_{\!\frac{i-j}{n-p}}\rpa\, \times\; \prod_{j=1}^{p-1}\; \B_{\frac{q_j}{n},\frac{j}{p} - \frac{q_j}{n}}\; \elaw\; \frac{n^n}{p^p}\; \prod_{j=0}^{p-1} \lpa \prod_{i=q_j+1}^{q_{j+1}-1} \G_{\frac{i}{n}}\rpa\, \times\; \prod_{j=1}^{p-1}\; \B_{\frac{q_j}{n},\frac{j}{p} - \frac{q_j}{n}},$$
where the second identity follows from
\begin{equation}
\label{BGA}
\B_{a, c-a}\times \G_c\, \elaw\, \G_a
\end{equation}
for all $c> a>0,$ which is a well-known consequence of (\ref{SF}). This completes the proof of the first claim, and it remains to show (\ref{BetaGamma}). To achieve this we use (2.6.20) in \cite{Z1} and again (\ref{LG}), which yields
$$\EE[{\bf Z}_{\frac{p}{n}}^{-ps}]\; =\;\frac{\Ga(1+ns)}{\Ga(1+ps)}\; =\; \lpa\frac{n^n}{p^p}\rpa^s \times\; \prod_{i=1}^{n-1}\frac{\Ga(\frac{i}{n} + s)}{\Ga(\frac{i}{n})}\;\times\;\prod_{j=1}^{p-1}\frac{\Ga(\frac{j}{p})}{\Ga(\frac{j}{p}+s)}$$
for every $s > -1/n.$ This can be rewritten
$$\EE[{\bf Z}_{\frac{p}{n}}^{-ps}]\; =\;\lpa\frac{n^n}{p^p}\rpa^s \times\; \prod_{j=0}^{p-1} \lpa \prod_{i=q_j+1}^{q_{j+1}-1} \frac{\Ga(\frac{i}{n} + s)}{\Ga(\frac{i}{n})}\rpa \times\; \prod_{j=1}^{p-1}\; \frac{\Ga(\frac{q_j}{n} + s)\Ga(\frac{j}{p})}{\Ga(\frac{j}{p} + s)\Ga(\frac{q_j}{n})},$$
and a fractional-moment identification based on (\ref{SF}) completes the proof of (\ref{BetaGamma}).

\endproof

\begin{REMS} {\em (a) When $p=1$ one recovers the well-known identity
\begin{equation}
\label{Will}
\Z_{\frac{1}{n}}^{-1}\;\elaw\; n^n\;\G_{\frac{1}{n}}\times\G_{\frac{2}{n}}\times\cdots\times\G_{\frac{n-1}{n}}
\end{equation}
for every $n\ge 2,$ which has been known since Williams \cite{W}. \smallskip

(b) When $p> 1$, other factorizations than (\ref{BetaGamma}) with $(p-1)$ Beta and $(n-p)$ Gamma random variables are possible, in choosing different indices in $[1,n-1]$ from the above $q_j, \, 1\le j\le p-1.$ This leads to a different localization of the Beta random variables inside the product - see Lemma 2 in \cite{TS1}. The above  choice was made in order to have Gamma random variables with parameters as small as possible, which will be important in the sequel.\smallskip

(c) In view of the identity $\L\elaw -\log [\B_{1,1}],$ the factorization (\ref{Beta}) is actually expressed in terms of Beta random variables only, and will hence be referred to as the "Beta factorization" subsequently. Contrary to (\ref{BetaGamma}), this factorization is canonical. For example, it leads after some simple rearrangements to the companion identities 
$$\Z_{\frac{p}{n}}^{-p}\;\elaw\; \L^{n-p}\; \times\; \K_{n,p}\qquad\mbox{and}\qquad \Z_{\frac{n-p}{n}}^{-(n-p)}\;\elaw\; \L^p\; \times\; \K_{n,p},$$
where 
$$\K_{n,p}\; \elaw\;\frac{n^n}{p^p (n-p)^{n-p}} \prod_{j=0}^{p-1} \lpa \prod_{i=q_j+1}^{q_{j+1}-1} \B_{\frac{i}{n},\frac{i-j}{n-p} - \frac{i}{n}}\rpa\; \times\; \prod_{j=1}^{p-1}\; \B_{\frac{q_j}{n},\frac{j}{p} - \frac{q_j}{n}}.$$
}
\end{REMS}

\section{Some properties of Kanter's random variable}

Kanter - see Corollary 4.1 in \cite{K} - observed the following independent factorization of $\Za$, which we express in terms of the Mittag-Leffler random variable: for every $\a\in (0,1)$ one has
\begin{equation}
\label{Kant}
\M_\a\; \elaw\; \L^{1-\a}\,\times\; b_\a(\U),
\end{equation}
where, here and throughout, $\U$ is the uniform random variable on $(0,1)$ and 
$$b_\a (u)\; = \;\frac{\sin (\pi u)}{\sin^\a (\pi\a u)\sin^{1-\a} (\pi(1-\a) u)}$$
for all $u\in (0,1).$ In the following, we will denote 
$$\Ka\; = \;b_\a(\U)$$ 
by the "Kanter random variable". It is interesting to mention in passing that the latter appears in the distributional theory of free stable laws - see the second part of Proposition A1.4 in \cite{BP} p.1054, and also \cite{ND} for related results. 

Notice that because $b_\a$ decreases from $\a^{-\a}(1-\a)^{\a-1}$ to $0$ (see the proof of Theorem 4.1 in \cite{K} for this latter fact) $\Ka$ is a bounded random variable with support $[0,\a^{-\a}(1-\a)^{\a-1}].$ In this section we first describe some further  distributional properties of $\Ka,$ which have their own interest and will partly play some r\^ole in the sequel. In the second paragraph, we prove some stochastic and convex orderings.

\subsection{Distributional properties}

We begin with the density function of $\Ka.$

\begin{PROP} 
\label{DK}
The density function of $\Ka$ is increasing and maps $(0, \a^{-\a}(1-\a)^{\a-1})$ onto $(1/(\Ga(\a)\Ga(1-\a)), +\infty)$.
\end{PROP}

\proof It follows from the proof of Lemma 2.1 in \cite{TS3} that $b_\a$ is strictly concave on $(0,1),$ whence the increasing character of the density of $\Ka.$ Besides one can compute
$$b_\a'(1-)\; =\; \frac{-\pi}{\sin(\pi\a)}\; = -\Ga(\a)\Ga(1-\a)\qquad \mbox{and}\qquad b_\a'(0+)\; =\; 0,$$
which entails that this density maps $(0, \a^{-\a}(1-\a)^{\a-1})$ onto $(1/(\Ga(\a)\Ga(1-\a)), +\infty).$

\endproof

\begin{REMS} {\em (a)  A further computation yields
$$b_\a''(1-)\; =\; \frac{2\pi(2\a -1)\cos(\pi\a)}{\sin^2(\pi\a)},$$
which is negative if $\a \neq 1/2$ and vanishes if $\a = 1/2.$ This shows that the derivative at zero of the density function of $\K_\a$ is positive if $\a \neq 1/2$ and vanishes if $\a = 1/2.$\smallskip

(b) Computing
$$b_\a''(0+)\; =\; \frac{\pi^2}{3\a^\a(1-\a)^{1-\a}}(2+3\a(1-\a)),$$
we observe that is not always smaller than $b_\a''(1-),$ so that $b_\a'$ is not convex in general. The latter would have entailed that $\Ka$ has a strictly convex density, a property which we believe to hold true notwithstanding.}
\end{REMS}

We next establish some identities in law, connecting $\Ka$ with Beta distributions and more generally with certain random variables characterized by their binomial moments, recently introduced in \cite{MP}. More precisely, it is shown in \cite{MP} that the sequence
$$\lpa\!\!\begin{array}{c} np + r \\ n\end{array}\!\!\rpa\; =\; \frac{\Ga(1+np +r)}{\Ga(1+n(p-1) +r)n!}$$
is positive definite for $p\ge 1$ and $-1\le r \le p-1$ and that it corresponds to the entire non-negative moments of some bounded random variable $X_{p,r}.$ The following proposition identifies the variable $X_{p,0}$ for all $p > 1.$

\begin{PROP} 
\label{3ID}
With the above notation, for any $\a\in (0,1)$ one has
\begin{equation}
\label{KMP}
\Ka\; \elaw\; {\bf K}_{1-\a}\; \elaw\; X_{1/\a, 0}^\a \,.
\end{equation}
Furthermore, when $\a = p/n$ is rational, then 
\begin{equation}
\label{KBye}
{\bf K}_{\frac{p}{n}}\;\elaw\;\lpa\frac{n^n}{p^p (n-p)^{n-p}} \prod_{j=0}^{p-1} \lpa \prod_{i=q_j+1}^{q_{j+1}-1} \B_{\frac{i}{n},\frac{i-j}{n-p} - \frac{i}{n}}\rpa\; \times\; \prod_{j=1}^{p-1}\; \B_{\frac{q_j}{n},\frac{j}{p} - \frac{q_j}{n}}\rpa^{\frac{1}{n}}
\end{equation}
with the above notation.
\end{PROP}

\proof The first identity in (\ref{KMP}) is plain because $b_\a \equiv b_{1-\a}.$ To obtain the second one, it is enough to compute the entire positive moments of $\Ka^{1/\a},$ and (\ref{Kant}) yields 
$$\EE[\Ka^{n/\a}]\; =\; \frac{\EE[\Za^{-n}]}{\EE[\L^{n(1-\a)/\a}]}\; =\; \frac{\Ga(1+n/\a)}{\Ga(1+n(1/\a-1))n!}\; =\;\EE[X_{1/\a,0}^n]$$
for every $n\ge 1.$ Last, the identity (\ref{KBye}) follows at once in comparing (\ref{Beta}) and (\ref{Kant}).

\endproof

\begin{REMS} {\em (a) The first identity in (\ref{KMP}) provides a direct explanation of  the above Remark 1 (c). The second identity in (\ref{KMP}) extends to $\a = 1$, where both sides are the deterministic random variable $\Un.$ \smallskip

(b) The identity (\ref{KMP}) shows that $X_{p,0} \elaw f_p(\U)$ for any $p\ge 1,$ where
$$f_p(u)\; =\; \frac{\sin^p (\pi u)}{\sin (\pi u/p) \sin^{p-1} ((1-1/p)(\pi u))}\cdot$$
It would be interesting to know whether such explicit representations exist for $X_{p,r}$ with $r\neq 0.$ Observe that Proposition 6.3 in \cite{MP} also represents $X_{p,0}$ as a free convolution power of the Bernoulli distribution.\smallskip

(c) Comparing (\ref{KMP}) and (\ref{KBye}) yields 
$$X_{\frac{n}{p}, 0}\;\elaw\;\lpa\frac{n^n}{p^p (n-p)^{n-p}} \prod_{j=0}^{p-1} \lpa \prod_{i=q_j+1}^{q_{j+1}-1} \B_{\frac{i}{n},\frac{i-j}{n-p} - \frac{i}{n}}\rpa\; \times\; \prod_{j=1}^{p-1}\; \B_{\frac{q_j}{n},\frac{j}{p} - \frac{q_j}{n}}\rpa^{\frac{1}{p}}$$
for every integers $n > p\ge 1.$ This is basically Theorem 3.3 in \cite{MP}, in the case $r=0.$}

\end{REMS}

We conclude with an identity in law which we will use during the proof of Theorem C.

\begin{PROP} 
\label{Joe}
For every $0<\bt<\a < 1$ one has
$${\bf Z}_{\frac{1-\a}{1-\bt}}^{\a -1}\,\times\, {\bf K}_\bt\; \elaw\; {\bf Z}_{\frac{\bt}{\a}}^{-\bt}\,\times\, {\bf K}_\a.$$
\end{PROP}

\proof This follows in comparing the fractional moments, which are given by
$$\frac{\Ga(1+s)}{\Ga(1+(1-\a)s)\Ga(1+\bt s)}$$
on both sides.

\endproof

\subsection{Comparison properties} In this paragraph we show that $\{\Ka^{\frac{1}{1-\a}},\, \a\in[0,1)\}$ and $\{\Ka,\, \a\in[0,1/2]\}$ can be arranged for the stochastic and convex orders, after suitable normalizations. This is the main technical contribution of the paper.

\begin{THE} 
\label{stcxK}
For every $0\le\bt<\a \le 1/2$ one has
\begin{equation}
\label{stK}
\a^\a(1-\a)^{1-\a} {\bf K}_\a\;\prost\; \bt^\bt(1-\bt)^{1-\bt}{\bf K}_\bt
\end{equation}
and
\begin{equation}
\label{cxK}
\Ga(1+\bt)\Ga(2-\bt) {\bf K}_\bt\;\prcvx\; \Ga(1+\a)\Ga(2-\a) {\bf K}_\a.
\end{equation}
For every $0\le\bt<\a <1$ one has
\begin{equation}
\label{stKa}
\a^{\frac{\a}{1-\a}}(1-\a) {\bf K}_\a^{\frac{1}{1-\a}}\;\prost\; \bt^{\frac{\bt}{1-\bt}}(1-\bt){\bf K}_\bt^{\frac{1}{1-\bt}}
\end{equation}
and
\begin{equation}
\label{cxKa}
(1-\bt) {\bf K}_\bt^{\frac{1}{1-\bt}}\;\prcvx\; (1-\a){\bf K}_\a^{\frac{1}{1-\a}}.
\end{equation}
\end{THE}

\begin{REMS}{\em In (\ref{stK}) and (\ref{cxK}) the requirement that $\a, \bt\in[0,1/2]$ is not a restriction, since by the first identity in (\ref{KMP}) all involved random variables have a parametrization which is symmetric w.r.t. $1/2.$}
\end{REMS}

To prove the theorem, we need the following lemmas.

\begin{LEM}
\label{lcxx}
For every $0\le\bt<\a \le 1/2,$ the function
$$x\;\mapsto\; \frac{\sin^\a (\pi\a x)\sin^{1-\a} (\pi(1-\a) x)}{\sin^{\bt} (\pi\bt x)\sin^{1-\bt} (\pi(1-\bt) x)}$$
is strictly log-convex on $(0,1).$ For every $0\le\bt<\a < 1,$ the function
$$x\;\mapsto\; \frac{\sin^{\frac{1}{1-\bt}} (\pi x)\sin^{\frac{\a}{1-\a}} (\pi\a x)\sin (\pi(1-\a) x)}{\sin^{\frac{1}{1-\a}} (\pi x)\sin^{\frac{\bt}{1-\bt}} (\pi\bt x)\sin (\pi(1-\bt) x)}$$
is strictly log-convex on $(0,1).$
\end{LEM}

\proof We begin with the first function, whose second logarithmic derivative equals
$$\pi^2\lpa\frac{\bt^3}{\sin^2 (\pi\bt x)}\, +\, \frac{(1-\bt)^3}{\sin^2 (\pi(1-\bt) x)}\, -\, \frac{\a^3}{\sin^2 (\pi\a x)}\, -\, \frac{(1-\a)^3}{\sin^2 (\pi(1-\a) x)}\rpa.$$
This can be rearranged into the sum of
$$\pi^2\lpa(1-\bt)\lpa\frac{(1-\bt)^2}{\sin^2 (\pi(1-\bt) x)}\, -\, \frac{\a^2}{\sin^2 (\pi\a x)}\rpa\, -\, (1-\a)\lpa\frac{(1-\a)^2}{\sin^2 (\pi(1-\a) x)}\, -\, \frac{\bt^2}{\sin^2 (\pi\bt x)}\rpa\rpa$$
and
$$\pi^2(1-\a-\bt)\lpa\frac{\a^2}{\sin^2 (\pi\a x)}\, -\, \frac{\bt^2}{\sin^2 (\pi\bt x)}\rpa.$$
For every $x\in(0,1)$ it can be checked that  on $(0,1)$ the function
\begin{equation}
\label{fcx}
t\;\mapsto\;\frac{t^2}{\sin^2 (\pi t x)}
\end{equation}
increases, which yields the positivity of the second summand since $\a+\bt< 1$, and is convex, whence the positivity of the first summand because $1-\bt> 1-\a.$ 

The argument for the second function is analogous. After some rearrangements, its second logarithmic derivative decomposes into the sum of
$$\pi^2\lpa\frac{1}{1-\a}\lpa\frac{1}{\sin^2 (\pi x)}\, -\, \frac{\a^2}{\sin^2 (\pi\a x)}\rpa\, -\, \frac{1}{1-\bt}\lpa\frac{1}{\sin^2 (\pi x)}\, -\, \frac{\bt^2}{\sin^2 (\pi\bt x)}\rpa\rpa$$
and
$$\pi^2\lpa\lpa\frac{(1-\bt)^2}{\sin^2 (\pi(1-\bt) x)}\, -\, \frac{(1-\a)^2}{\sin^2 (\pi(1-\a) x)}\rpa\, +\, \lpa\frac{\a^2}{\sin^2 (\pi\a x)}\, -\, \frac{\bt^2}{\sin^2 (\pi\bt x)}\rpa\rpa.$$
Again, the positivity of the first summand comes from the convexity of the function (\ref{fcx}), whereas the positivity of the second summand follows from its increasing character.

\endproof

\begin{LEM} 
\label{bound}
The functions 
$$t\;\mapsto\; \frac{t(1-t)}{\sin^2(\pi t)}\qquad \mbox{\em and}\qquad t\;\mapsto\;\frac{\sin(\pi t)}{t^{1-t}(1-t)^t}$$
decrease on $(0,1/2].$ The function $t\mapsto t^{\frac{t}{1-t}}$ decreases on $(0,1).$
\end{LEM}

\proof Let us start with the first function, whose logarithmic derivative is
$$\frac{1}{t}\, -\, \frac{1}{1-t}\, -\, 2\pi \cot(\pi t)$$
and vanishes on $t = 1/2.$ The second logarithmic derivative is
$$\frac{2\pi^2}{\sin^2(\pi t)}\, -\,\frac{1}{t^2}\, -\, \frac{1}{(1-t)^2}\; >\; \frac{2}{t^2}\lpa \frac{\pi^2t^2}{\sin^2(\pi t)}\, -\, 1 \rpa\; > \; 0$$
for all $t\in (0,1/2),$ which gives the claim. The argument for the second function is more involved. Its logarithmic derivative is
$$\frac{1}{1-t}\, -\, \frac{1}{t}\, +\,\log(t)\, -\, \log (1-t)\, +\, \pi \cot(\pi t)$$
and, again, vanishes on $t = 1/2.$ The second logarithmic derivative is
\begin{eqnarray*}
\frac{1}{t^2}\, +\, \frac{1}{(1-t)^2}\, +\, \frac{1}{t(1-t)} \,-\,\frac{\pi^2}{\sin^2(\pi t)} & = & \frac{1}{t(1-t)}\, -\, \sum_{n\ge 1} \frac{1}{(n+t)^2} \, -\,  \sum_{n\ge 2} \frac{1}{(n-t)^2} \\
& > & 4\, -\, \sum_{n\ge 1} \frac{1}{n^2} \, -\,  \sum_{n\ge 2} \frac{4}{(2n-1)^2}  \; =\; 8 \, -\, \frac{2\pi^2}{3} \; > \; 0,
\end{eqnarray*}
which finishes the claim. The increasing character of the third function is easy and we leave the details to the reader.

\endproof

\noindent
{\bf Proof of Theorem \ref{stcxK}.} Let us start with (\ref{stK}). Consider the function
$$x\; \mapsto\; \frac{\bt^\bt(1-\bt)^{1-\bt}b_\bt(x)}{\a^\a(1-\a)^{1-\a}b_\a(x)}$$
which is strictly log-convex on $(0,1)$ if $0\le\bt<\a \le 1/2,$ by Lemma \ref{lcxx}. At $0+$ its limit is 1, whereas the limit of its derivative is 0. Putting everything together  entails that this function is strictly greater than 1 on $(0,1).$ From the definition of $\Ka$ we deduce 
$$\a(1-\a)^{1-\a}\K_\a\; \prst\; \bt(1-\bt)^{1-\bt}\K_\bt,$$
whence (\ref{stK}) since it is also clear by construction that the constants are optimal. 

Let us now consider (\ref{cxK}), introducing the strictly log-convex function
$$x\; \mapsto\; \frac{\Ga(1+\bt)\Ga(2-\bt)b_\bt(x)}{\Ga(1+\a)\Ga(2-\a) b_\a(x)},$$
whose limits at $0+$
resp. $1-$ are
$$\frac{\bt^{1-\bt}(1-\bt)^\bt\sin(\pi \a)}{\a^{1-\a}(1-\a)^\a\sin(\pi \bt)} \; <\; 1\qquad \mbox{resp.}\qquad \frac{\bt(1-\bt)\sin^2(\pi \a)}{\a(1-\a)\sin^2(\pi \bt)} \; >\; 1,$$
by Lemma \ref{bound}. By convexity, we see that the distribution function of  $\Ga(1+\a)\Ga(2-\a){\bf K}_\a$ crosses that of $\Ga(1+\bt)\Ga(2-\bt){\bf K}_\bt$ at exactly one point, and from above. Since 
$$ \EE[\Ga(1+\a)\Ga(2-\a){\bf K}_\a]\, =\,  \EE[\Ga(1+\bt)\Ga(2-\bt){\bf K}_\bt]\, =\,1$$ 
for all $\a,\bt,$ this completes the proof of (\ref{cxK}) by Theorem 2.A.17 in \cite{SS}.

The arguments for (\ref{stKa}) and (\ref{cxKa}) are the same. To get (\ref{stKa}), consider the function 
$$x\; \mapsto\; \frac{(1-\bt)\bt^{\frac{\bt}{1-\bt}}b_\bt^{\frac{1}{1-\bt}}(x)}{(1-\a)\a^{\frac{\a}{1-\a}}b_\a^{\frac{1}{1-\a}}(x)},$$
which is strictly log-convex on $(0,1)$ if $0\le\bt<\a < 1$ by Lemma \ref{lcxx}, has limit 1 at $0+$, whereas the limit of its derivative is 0. To obtain (\ref{cxKa}), consider first the strictly log-convex function 
$$x\; \mapsto\; \frac{(1-\bt)b_\bt^{\frac{1}{1-\bt}}(x)}{(1-\a)b_\a^{\frac{1}{1-\a}}(x)},$$
whose limit at $0+$ is
$$\a^{\frac{\a}{1-\a}}\bt^{\frac{\bt}{\bt-1}}\; < \; 1$$
by Lemma \ref{bound}, whereas its limit at $1-$ is $+\infty.$ To finish the proof observe that
$$ \EE[(1-\a){\bf K}_\a^{\frac{1}{1-\a}}]\, =\,  \EE[(1-\bt){\bf K}_\bt^{\frac{1}{1-\bt}}]\, =\,1,$$
and use again Theorem 2.A.17 in \cite{SS}.

\qed

\section{Proof of Theorems A and B}

\subsection{Proof of Theorem A} It is plain from (\ref{Kant}) and (\ref{stKa}) that 
$$(1-\a)\a^{\frac{\a}{1-\a}} {\bf Z}_\a^{\frac{-\a}{1-\a}}\;\prst\;(1-\bt)\bt^{\frac{\bt}{1-\bt}}{\bf Z}_\bt^{\frac{-\bt}{1-\bt}}$$
for all $0<\bt<\a<1.$ On the other hand, it is well-known \cite{C} and easy to see that
\begin{equation}
\label{Cresson}
{\bf Z}_\bt^{-\bt}\;\claw\; \L\qquad\mbox{as $\bt \to 0.$}
\end{equation}
An application of Stirling's formula and (2.6.20) in \cite{Z1} yields
$$\EE\lcr(1-\a)^s {\bf Z}_\a^{\frac{-\a s}{1-\a}}\rcr\;\rightarrow\; s^s\qquad\mbox{as $\a \to 1$}$$
for every $s > 0,$ whence $(1-\a)\a^{\frac{\a}{1-\a}} {\bf Z}_\a^{\frac{-\a}{1-\a}} \claw \S$ as $\a \to 1.$ Putting everything together entails
$$\S \;\prst\;(1-\a)\a^{\frac{\a}{1-\a}} {\bf Z}_\a^{\frac{-\a}{1-\a}}\;\prst\;(1-\bt)\bt^{\frac{\bt}{1-\bt}}{\bf Z}_\bt^{\frac{-\bt}{1-\bt}}\;\prst\;\L$$
for all $0<\bt<\a<1.$ To show (\ref{stZ}), it remains to prove that there does not exist any $c > 1$ such that
$$\pb\lcr c\S \ge x\rcr\;\le\;\pb\lcr \L \ge x\rcr, \qquad x\ge 0.$$
But if the latter were true, then for any $s>0$ we would plainly have
$$\lpa \frac{cs}{e}\rpa^{s}\; =\; c^s\EE\lcr\S^s\rcr\;\le\;\EE\lcr \L^s \rcr\; =\; \Ga(1+s)$$
and this would contradict Stirling's formula. This completes the proof of (\ref{stZ}). The proof of (\ref{cxZ}) is a simple consequence of (\ref{Kant}), (\ref{cxKa}), and Theorem 2.A.6.(b) in \cite{SS}.

\qed

\begin{REMS} {\em Another way to see that the normalization is optimal in (\ref{stZ}) is to use the behaviour of the distribution function of $\Za$ at zero - see e.g. (14.31) p. 88 in \cite{S}:
\begin{equation}
\label{zero}
x^{\a/(1-\a)}\log \pb [\Za \le x] \to (1-\a)\a^{\a/(1-\a)} \quad \mbox{as $x\to 0.$}
\end{equation}}
\end{REMS}

\subsection{Proof of Theorem B} A repeated use of Theorem 2.A.6.(b) in \cite{SS} shows that it is enough to show that 
\begin{equation}
\label{cxL}
\Ga(2-\bt)\L^{1-\a}\;\prcvx\;\Ga(2-\a)\L^{1-\bt}
\end{equation}
for all $0\le \bt < \a \le 1.$ Indeed, by (\ref{Kant}) we have 
$$\M_u\; \elaw\; \frac{{\bf L}^{1-u}}{\Ga(2-u)}\, \times\, 
\Ga(1+u)\Ga(2-u) {\bf K}_u$$
for all $u\in [0,1]$ and we know from Theorem \ref{stcxK} and the first identity in law of Proposition \ref{3ID} that
$$\Ga(1+\a)\Ga(2-\a) {\bf K}_\a\;\prcvx\; \Ga(1+\bt)\Ga(2-\bt) {\bf K}_\bt$$
for all $1/2\le \bt < \a \le 1.$ To show (\ref{cxL}), observe that
$$\EE\lcr \Ga(2-\a)\L^{1-\bt}\rcr\; = \; \EE\lcr \Ga(2-\bt)\L^{1-\a}\rcr\; =\; \Ga(2-\a)\Ga(2-\bt)$$
and that 
$$\pb\lcr \Ga(2-\a)\L^{1-\bt}\le x\rcr\; \sim \; \lpa\frac{x}{\Ga(2-\a)}\rpa^{\frac{1}{1-\bt}}\; \gg \; \lpa\frac{x}{\Ga(2-\bt)}\rpa^{\frac{1}{1-\a}}\;\sim \;\pb\lcr \Gamma(2-\bt)\L^{1-\a}\le x\rcr$$
at zero. Therefore, again by Theorem 2.A.17 in \cite{SS} we need to show that the densities of the random variables $\Ga(2-\bt)\L^{1-\a}$ and $\Ga(2-\a)\L^{1-\bt}$ only meet twice. This amounts to $\sharp\{t >0, \; ct\, e^{-c t} = \gamma t^\gamma\, e^{-t^\gamma}\} = 2,$ where 
$$\gamma\; =\; \frac{1-\a}{1-\bt}\; \in \; (0,1)$$
and $c >0$ is some normalizing constant. The latter claim follows easily from the strict concavity of $t\mapsto(1-\gamma)\log(t)\, -\, c t \, +\, t^\gamma.$

\qed

\begin{REMS} {\em Consider the positive random variable $\I_\a = -\inf\{X_t, \; t\le 1\},$  where $\{X_t, \; t\ge 0\}$ is the spectrally negative strictly $(1/\a)-$stable L\'evy process introduced before the statement of Theorem B. Formula (9) in \cite{TS0} shows that 
$$\I_\a\;\elaw\;\M_\a \; \times\; \Y_\a,$$
where $\Y_\a$ has density
$$\frac{-\sin(\pi/\a) x^{\frac{1}{\a} -2}(1+x)}{\pi(x^{\frac{2}{\a}} - 2x^{\frac{1}{\a}}\cos(\pi/\a) +1)}$$
over $(0,+\infty).$ An application of Theorem 1 in \cite{TS0} entails $\EE[\Y_\a] =1$ for every $\a\in[1/2, 1),$ so that 
$$\EE[\Gamma(1+\a)\I_\a]\; =\;\EE[\Gamma(1+\a)\M_\a]  \;=\; 1.$$ 
It is hence natural to ask whether the family $\{\Gamma(1+\a)\I_\a, \;1/2\le \a < 1\}$ could also be arranged along the convex order. An analysis similar to the above shows that  $\Y_\bt\,\prcvx\,\Y_\a$ for every $\bt < \a,$ so that we cannot conclude directly. We believe that the ordering for $\Gamma(1+\a)\I_\a$ is in the same direction, that is
$$\Gamma(1+\bt)\I_\bt\;\prcvx\;\Gamma(1+\a)\I_\a$$ 
for every $1/2\le \bt < \a<1.$ This would allow to extend the martingale introduced after the statement of Theorem B. Observe that this martingale would not converge since $\I_\a$ does not have a limit law when $\a\to 1.$}
\end{REMS}

\section{Proof of Theorem C}

\subsection{The case $\a\le 1/2$} We need the following lemma.

\begin{LEM}
\label{BG}
For every integers $0<p<n,$ one has
\begin{equation}
\label{BGbis}
\B_{\frac{1}{n}, \frac{1}{p}-\frac{1}{n}}^{1/n}\; \prost\; p^{p/n}\, \Ga(1-p/n)\, \G_{\frac{1}{n}}^{1/n}\times\;\lpa\prod_{i=2}^p\, \G_{\frac{i}{p}}\rpa^{1/n}
\end{equation}
and
\begin{equation}
\label{BGter}
 \B_{\frac{1}{n}, \frac{1}{p}-\frac{1}{n}}^{1/n}\; \prcvx\; p^{p/n}\,\Ga(1+p/n)^{-1}\, \G_{\frac{1}{n}}^{1/n}\times\;\lpa\prod_{i=2}^p\, \G_{\frac{i}{p}}\rpa^{1/n}.
\end{equation}
\end{LEM}

\proof The density of the random variable on the left-hand side of (\ref{BGbis}) equals
$$\frac{n\Gamma(1/p)\, (1-x^n)^{1/p-1/n-1}}{\Gamma(1/n)\Gamma(1/p-1/n)} \,\Un_{(0,1)}(x)\,\to\, \frac{\Gamma(1/p)}{\Gamma(1/n+1)\Gamma(1/p-1/n)} \quad\mbox{as $x\to 0+,$}$$
and increases on $(0,1)$ because $1/p-1/n-1 <0.$ Setting $f_{n,p}$ for the density of
$$p^{p/n}\, \Ga(1-p/n)\,\times\lpa\prod_{i=2}^p \G_{\frac{i}{p}}\rpa^{1/n}$$
we see by multiplicative convolution that the density of the random variable on the right-hand side of (\ref{BGbis}) equals
$$\frac{n}{\Gamma(1/n)} \int_0^\infty e^{-(x^n/y)} \, f_{n,p}(y) \frac{dy}{y}$$
and is hence positive decreasing on $(0,+\infty).$  Besides, its value at $0+$ is
\begin{eqnarray*}
\frac{n}{p^{p/n}\, \Gamma(1/n)\Ga(1-p/n)}\,\times\EE\lcr\prod_{i=2}^p \G_{\frac{i}{p}}^{-1/n}\rcr, 
\end{eqnarray*}
which transforms by (\ref{LG}) into 
$$\frac{n\Gamma(1/p)}{\Gamma(1/n)\Gamma(1/p-1/n)} \, \times\, \frac{\Gamma(1/p-1/n)\times\cdots\times\Gamma(1-1/n)}{p^{p/n}\Gamma(1/p-1/n)\Gamma(1/p)\times\cdots\times\Gamma(1)}\; =\;  \frac{n\Gamma(1/p)}{\Gamma(1/n)\Gamma(1/p-1/n)}\cdot$$
Putting everything together and plotting the two densities, we easily deduce (\ref{BGbis}).\smallskip

To obtain (\ref{BGter}), observe that the density of the random variable on the left-hand side is positive decreasing on $(0,+\infty)$ and equals
\begin{eqnarray*}
\Gamma(1+p/n)\Gamma(1-p/n)\,\times\,\frac{n\Gamma(1/p)}{\Gamma(1/n)\Gamma(1/p-1/n)} &=& \frac{\pi p\Gamma(1/p)}{\sin(\pi p/n)\Gamma(1/n)\Gamma(1/p-1/n)}\\ 
&> & \frac{n\Gamma(1/p)}{\Gamma(1/n)\Gamma(1/p-1/n)}
\end{eqnarray*}
at $0+,$ by Euler's reflection formula for the Gamma function. Since the increasing density of the random variable on the right-hand side tends to $+\infty$ at $1-$ and equals
$$\frac{n\Gamma(1/p)}{\Gamma(1/n)\Gamma(1/p-1/n)}$$
at $0+$, this means that the two densities meet at exactly one point on $(0,1).$ Computing with the help of (\ref{LG}) 
\begin{eqnarray*}
\EE\lcr p^{p/n}\, \Ga(1+p/n)^{-1}\, \G_{\frac{1}{n}}^{1/n}\times\;\lpa\prod_{i=2}^p\, \G_{\frac{i}{p}}\rpa^{1/n}\rcr & = & \frac{\Ga(2/n)p^{p/n}\Ga(2/p+1/n)\cdots \Ga(p/p+1/n)}{\Ga(1/n)\Ga(1+p/n)\Ga(2/p)\cdots \Ga(p/p)}\\
& = & \frac{\Ga(2/n)\Ga(1/p)}{\Ga(1/n)\Ga(1/p+1/n)}\; =\; \EE\lcr\B_{\frac{1}{n}, \frac{1}{p}-\frac{1}{n}}^{1/n}\rcr
\end{eqnarray*}
and plotting the two densities, we finally get (\ref{BGter}) from Theorem 2.A.17 in \cite{SS}.

\endproof

\noindent
{\bf End of the proof.} Let us begin with the stochastic order. By Formula (14.35) p.88 in \cite{S} we have
\begin{equation}
\label{ZAL}
\pb[\M_\a\le x] \;=\; \pb[\Za \ge x^{-1/\a}]\;\sim\; \frac{x}{\Ga(1-\a)}\; \sim\; \pb [\Ga(1-\a)\L\le x]
\end{equation}
as $x\to 0+,$ and it is hence enough to show that
$$\Z_\a\;\prst\; \Ga(1-\a)\L.$$
To do so, we first suppose that $\a=p/n$ for some positive integers $p,n$ with $n\ge 2p.$ For the sake of clarity we will consider the case $p=1$ separately. By (\ref{BGbis}) we then have
$$\B_{1/n, 1-1/n}\; \prst\; \Ga(1-1/n)^n\, \G_{1/n},$$
whence we deduce, by (\ref{BGA}) and Theorem 1.A.3 (d) in \cite{SS},
$$\G_{1/n}\; \prst\; \Ga(1-1/n)^n\, \G_{1/n}\times \G_1.$$
By (\ref{BGA}) and again Theorem 1.A.3 (d) in \cite{SS}, this yields
$$\Z_{1/n}^{-1}\, \elaw\, n^n \G_{1/n}\times\G_{2/n}\times \cdots\times \G_{(n-1)/n}\, \prst\, n^n\Ga(1-1/n)^n\, \G_{1/n}\times\G_{2/n}\times \cdots\times \G_1\,\elaw\, \Ga(1-1/n)^n \L^n,$$
where the identity in law follows from the first identity in (\ref{IDL}). We finally obtain the required claim
$$\M_{1/n}\; \prst\; \Ga(1-1/n) \L.$$
We now consider the case $p\ge 2.$ Since $n\ge 2p$ by assumption, observe that $q_1\ge 2$ with the notation of Section 2. By (\ref{BetaGamma}) and (\ref{BGA}), we can rewrite
\begin{eqnarray*}
{\bf Z}_{\frac{p}{n}}^{-p} & \elaw & \frac{n^n}{p^p}\; \G_{\frac{1}{n}}\,\times \lpa\prod_{i=2}^{q_2-1} \G_{\frac{i}{n}}\rpa\times \;\prod_{j=1}^{p-1} \lpa\lpa \prod_{i=q_j+1}^{q_{j+1}-1} \G_{\frac{i}{n}}\rpa\, \times\, \B_{\frac{q_j}{n},\frac{j}{p} - \frac{q_j}{n}}\rpa\\
& \elaw & \frac{n^n}{p^p}\; \B_{\frac{1}{n}, \frac{1}{p}-\frac{1}{n}}\,\times \, \G_{\frac{1}{p}} \,\times \lpa\prod_{i=2}^{q_2-1} \G_{\frac{i}{n}}\rpa\times \;\prod_{j=1}^{p-1} \lpa\lpa \prod_{i=q_j+1}^{q_{j+1}-1} \G_{\frac{i}{n}}\rpa\, \times\, \B_{\frac{q_j}{n},\frac{j}{p} - \frac{q_j}{n}}\rpa.
\end{eqnarray*}
On the other hand, it follows from a repeated use of (\ref{BGA}) and the first identity in (\ref{IDL}) that
\begin{eqnarray*}
\L^n & \elaw & n^n\,  \G_{\frac{1}{n}}\,\times \lpa\prod_{i=2}^{q_2-1} \G_{\frac{i}{n}}\rpa\times \;\prod_{j=1}^{p-1} \lpa\lpa \prod_{i=q_j+1}^{q_{j+1}-1} \G_{\frac{i}{n}}\rpa\, \times\, \G_{\frac{q_j}{n}}\rpa\; \times\; \L\\
& \elaw & n^n\,  \G_{\frac{1}{n}}\,\times \, \lpa\prod_{i=2}^p\, \G_{\frac{i}{p}}\rpa\,\times\;\G_{\frac{1}{p}}\,\times \lpa\prod_{i=2}^{q_2-1} \G_{\frac{i}{n}}\rpa\times \;\prod_{j=1}^{p-1} \lpa\lpa \prod_{i=q_j+1}^{q_{j+1}-1} \G_{\frac{i}{n}}\rpa\, \times\, \B_{\frac{q_j}{n},\frac{j}{p} - \frac{q_j}{n}}\rpa\cdot
\end{eqnarray*}
By (\ref{BGbis}) and Theorem 1.A.3 (d) in \cite{SS} we deduce
$$\M_{p/n}\; \prst\; \Ga(1-p/n) \L,$$
which is the required claim for $\a = p/n$ with $n\ge 2p.$ The general case for all $\a\in(0,1/2]$ follows by density. The proof of 
$$ \Ga(1+\a) \M_\a\; \prcvx\;\L$$
for all $\a\in(0,1/2]$ goes along the same lines in using (\ref{BGter}), Theorem 2.A.6 (b) in \cite{SS}, and a density argument. We leave the details to the reader.

\qed

\begin{REMS} {\em The same kind of proof can also be performed to handle the case $\a > 1/2$, with different and more complicated details.}
\end{REMS}

\subsection{The case $\alpha > 1/2$} We will use an argument closer to that of Theorems A and B, in order to get a more general result which is stated in Theorem \ref{Mike} below. 

\begin{LEM} 
\label{Pete}
With the above notation, the random variable
$$\Z_{\frac{\bt}{\a}}^{-\bt}\times\Ka$$
has a non-increasing density for every $0<\bt \le 1/2$ and $\a \ge 1/2\vee(2\bt\wedge (\bt+1)/2).$
\end{LEM}

\proof We first suppose that $0<\bt\le \a/2<1/2\le\a<1.$ The required property means that the distribution function of $\Z_{\bt/\a}^{-\bt}\times\Ka$ is convex and, by a density argument and because the pointwise limit of convex functions is convex, this allows to consider the case where $\a,\bt$ are rational. We hence set
$$\frac{\bt}{\a}\;=\;\frac{k}{l}\; \le \; \frac{1}{2}, \qquad \bt\;=\;\frac{p}{q}\; \le \; \frac{1}{2}, \qquad\mbox{and}\qquad \a\;=\;\frac{lp}{kq}\; \ge \; \frac{1}{2}\cdot$$
Since $k/l\le 1/2,$ by (\ref{BetaGamma}) we first observe that $\Z_{k/l}^{-p/q}$ admits 
$\G_{1/l}^{p/kq}$ as a multiplicative factor. Reasoning as for (\ref{IDL}), we obtain that the latter factorizes with  $\G_{1/lp}^{1/kq}.$ On the other hand, using (\ref{KBye}) and observing that $q_1 = 1$ therein because $lp\ge kq/2,$ we see that $\K_{lp/kq}$ has $(\B_{1/kq, 1/lp- 1/kq})^{1/kq}$ as a multiplicative factor. Hence, the random variable $\Z_{\bt/\a}^{-\bt}\times\Ka=\Z_{k/l}^{-p/q}\times\K_{lp/kq}$ factorizes with 
$$\lpa\G_{\frac{1}{lp}}\times\B_{\frac{1}{kq}, \frac{1}{lp}- \frac{1}{kq}}\rpa^{\frac{1}{kq}}\; \elaw\; \G_{\frac{1}{kq}}^{\frac{1}{kq}}$$
and has, reasoning as in Lemma \ref{BG}, a non-increasing density. To obtain the same property in the case $0<\bt\le 1/2<(\bt+1)/2\le\a,$ it suffices to apply Proposition \ref{Joe}. 

\endproof

\begin{REMS}{\em By Proposition \ref{DK} and the main theorem in \cite{TS1}, we know that the random variable 
$$\Z_{\frac{\bt}{\a}}^{-\bt}\times\Ka$$
is unimodal as soon as $2\bt\le \a.$ However, its mode is positive when $\bt > 1/2.$ Indeed, if its mode were zero, then the random variable
$$\L^{1-\a}\times\Z_{\frac{\bt}{\a}}^{-\bt}\times\Ka\;\elaw\; \lpa\Za\times\Z_{\frac{\bt}{\a}}^{\frac{\bt}{\a}} \rpa^{-\a}\; \elaw\; \M_\bt$$
would also have a non-increasing density, which is false - see again (14.35) p. 88 in \cite{S}.}
\end{REMS}

We can now state the main result of this paragraph, which finishes the proof of Theorem C in the case $\a > 1/2$ in letting $\bt\to 0$ and using (\ref{Cresson}).

\begin{THE}
\label{Mike}
For every $0<\bt \le 1/2$ and $\a \ge 1/2\vee(2\bt\wedge (\bt+1)/2)$
one has
\begin{equation}
\label{cxstZa}
\Gamma(1-\bt)\M_\a\;\prost\;\Gamma(1-\a)\M_\bt\qquad\mbox{and}\qquad \Gamma(1+\a)\M_\a\;\prcvx\; \Gamma(1+\bt)\M_\bt.
\end{equation}
\end{THE}

\proof Let us begin with the stochastic order. We have to compare
$$\Gamma(1-\bt)\M_\a\; \elaw\; \Gamma(1-\bt)\L^{1-\a}\,\times\,\Ka$$
and
\begin{eqnarray*}
\Gamma(1-\a)\M_\bt \; \elaw\; \Gamma(1-\a)\L^{1-\bt}\,\times\,\K_\bt & \elaw & \Gamma(1-\a)\Z_{\frac{1-\a}{1-\bt}}^{\a-1}\times \,\L^{1-\a}\,\times\,\K_\bt \\
& \elaw & \Gamma(1-\a)\L^{1-\a}\,\times\,\Z_{\frac{\bt}{\a}}^{-\bt}\times \,\Ka,
\end{eqnarray*}
where the second identity in law follows from Shanbhag-Sreehari's identity - see e.g. Exercise 29.16 in \cite{S} - and the third one from Proposition \ref{Joe}. Using again (\ref{ZAL}), we see that we are reduced to prove 
$$\Gamma(1-\bt)\Ka \; \prst\; \Gamma(1-\a)\Z_{\frac{\bt}{\a}}^{-\bt}\times \,\Ka.$$
By Proposition \ref{DK}, the density of the random variable on the left-hand side increases on $(0,\Gamma(1-\bt))$ whereas by Lemma \ref{Pete}, the density of the random variable on the right-hand side is non-increasing on $(0, +\infty).$ Hence, reasoning as above, we need to show that the two densities coincide at $0\!+\!.$ By Proposition \ref{DK}, the value for the left-hand side is
$$\frac{1}{\Ga(\a)\Ga(1-\a)\Ga(1-\bt)},$$
whereas again by (2.6.20) in \cite{Z1} the evaluation for the right-hand side is
$$\frac{1}{\Ga(\a)\Ga(1-\a)^2}\times \EE[\Z_{\frac{\bt}{\a}}^{\bt}]\; =\; \frac{1}{\Ga(\a)\Ga(1-\a)^2}\times \frac{\Ga(1-\a)}{\Ga(1-\bt)}\; =\; \frac{1}{\Ga(\a)\Ga(1-\a)\Ga(1-\bt)}\cdot$$
The proof of the convex ordering goes along the same lines as for (\ref{BGter}), using  
$$\EE[ \Gamma(1+\bt)\M_\bt]\;=\; \EE[\Gamma(1+\a)\M_\a]\; =\; 1,$$
and we leave the details to the reader.

\endproof

\begin{REMS} 
\label{F1}
{\em We believe that the two statements in (\ref{cxstZa}) are true without any restriction on $0< \bt< \a \le 1.$ Observe that Theorem B already shows the convex ordering on $1/2\le \bt< \a \le 1.$ Using (\ref{Will}) and an adaptation of Lemma \ref{BG}, it is possible to show that
$$\Gamma(1-1/n)\M_{\frac{1}{p}} \;\prost\;\Gamma(1-1/p)\M_{\frac{1}{n}}\quad\mbox{and}\quad \Gamma(1+1/p)\M_{\frac{1}{p}}\;\prcvx\; \Gamma(1+1/n)\M_{\frac{1}{n}}$$
for every integers $1< p<n$. Unfortunately, this method does not work to get the inequalities in general since the Beta factorization does not simplify enough for this purpose, and also because of Remark 7.}
\end{REMS}

\section{Applications of the main results}

\subsection{Behaviour of positive stable medians} If $X$ is a real random variable whose density does not vanish on its support, we denote its median by $m^{}_X,$ which is the unique number $m_X$ such that
$$\pb[X \le m^{}_X]\; =\; \pb[X \ge m^{}_X]\; =\; \frac{1}{2}\cdot$$
Notice that if $X\,\prst\, Y$, then $m^{}_X\,\le\, m^{}_Y.$ The number $m_\a = m^{}_{\Za}$ is not explicit except in the case $\a=1/2$ with
$$m_{1/2}\; =\; \frac{1}{4 ({\rm Erfc}^{-1}(1/2))^2}\; \sim\; 1.0990,$$ 
where ${\rm Erfc}^{-1}$ stands for the inverse of the complementary error function. Observe that (\ref{stZ}) entails 
$$0\; <\; m_\S\; \le \; \frac{1}{4 m_{1/2}}\; =\; ({\rm Erfc}^{-1}(1/2))^2\; \sim\; 0.2274,$$ which is equivalent to $m_{\X+1} = \log (m_\S)+1\le -0.481.$ This latter estimate on the median of the completely asymmetric Cauchy random variable $\X+1$ does not seem to follow directly from Formula (2.2.19) in \cite{Z1}. It is not even clear from this formula whether $m_{\X+1} < 0$ viz. $\pb[\X + 1 \ge 0] > 1/2,$ which is a statement on the positivity parameter. 

A  combination of (\ref{stZ}) and (\ref{cxstL}) and the fact that $m_\L =\log (2)$ yield the following bounds on $m_\a$ for every $\a\in(0,1):$
\begin{equation}
\label{Bound}
\a\lpa\frac{1-\a}{\log(2)}\rpa^{\frac{1-\a}{\a}}\vee\; \lpa \frac{1}{\log(2) \Ga(1-\a)}\rpa^{\frac{1}{\a}} \; \le \; m_\a\; \le\; \a\lpa \frac{1-\a}{m_\S}\rpa^{\frac{1-\a}{\a}}\cdot
\end{equation}
These bounds show that $m_\a\to +\infty$ at exponential speed when $\a\to 0\!+\!.$ On the other hand, (\ref{Bound}) also entails that $m_\a\le \a$ for every $\a\in(1-m_\S, 1)$ and since it is clear that $m_\a\to 1$ as $\a\to 1-,$ overall we observe the curious fact that the function $\a\to m_\a$ is not monotonic on $(0,1).$ This is in sharp contrast with the behaviour of the mode of $\Za$ which is, at least heuristically, an increasing function of $\a$ - see the final tables in \cite{NS}. The following proposition shows however a monotonic behaviour for $m_\a$ in the neighbourhood of $1.$ 

\begin{PROP}
\label{Muni}
With the above notation, the function $\a\mapsto m_\a$ increases on $(1-m_\S, 1).$ 
\end{PROP}

\proof Setting $f_\a(x) = (1-\a)\a^{\frac{\a}{1-\a}} x^{\frac{-\a}{1-\a}}$ for every $x > 0,$ it is clear from (\ref{stZ}) that 
$$f_\a(m_\a)\;\le\; f_\bt(m_\bt)$$
for every $0<\bt<\a<1.$ Because $f_\a$ decreases on $(0,+\infty)$, it is enough to show that $f_\bt(m_\bt) <  f_\a(m_\bt)$ for every $1-m_\S <\bt<\a<1.$ Since then we have $m_\bt\le \bt$ by the previous observation, this amounts to show that
$$t\; \mapsto\; f_t(x)$$
increases on $(x,1)$ for every $x\in (0,1).$ Computing the logarithmic derivative
$$\frac{\log(t) -\log(x)}{(1-t)^2}$$
finishes the proof.

\endproof

\begin{REMS} {\em It would be interesting to investigate the global minimum of the function $\a\mapsto m_\a.$  We believe that this function decreases and then increases.}
\end{REMS} 

Our next result shows a partial median-mode inequality for $\Za,$ which is actually a mean-median-mode inequality because $\EE[\Za] =+\infty.$ The reader can consult \cite{BDG} for more details and references concerning mean-median-mode inequalities. Setting $M_\a$ for the mode of $\Za,$ that is the unique local maximum of its density, we recall that $M_\a$ is explicit only for $\a = 1/2$ with $M_{1/2} =1/6.$ We refer to \cite{NS} and the citations therein for more information on $M_\a.$

\begin{PROP}
\label{2M}
With the above notation, one has
$$M_\a \; < \; m_\a$$
as soon as $\a < 1/(1+\log(2)) \sim 0.5906.$
\end{PROP}

\proof We use the following upper bound  
$$M_\a\; \le\; \lpa \frac{\a}{\Gamma(2-\a)}\rpa^{\frac{1}{\a}},$$
which is quoted in \cite{NS} as a consequence of (6.4) in \cite{SY}. Combining with (\ref{Bound}) we need to find the set of $\a\in(0,1)$ such that $\a\log(2) \Ga(1-\a)\, <\,\Gamma(2-\a).$ This set is $(0,1/(1+\log(2)).$

\endproof

\begin{REMS}{\em (a) The first lower bound in (\ref{Bound}) behaves better than the second one when $\a\to 1,$ and can be improved by (\ref{stZ}) into 
\begin{equation}
\label{Best}
m_\a \; \ge \; \a(4m_{1/2}(1-\a))^{\frac{1-\a}{\a}}
\end{equation}
for every $\a\in (1/2,1).$ Unfortunately, this does not extend the validity domain in the above proposition. Indeed, it follows from the log-convexity of the Gamma function that
$$\lpa\frac{\a}{\Gamma(2-\a)}\rpa^{\frac{1}{\a}}\;>\;\a(4m_{1/2}(1-\a))^{\frac{1-\a}{\a}}$$
for every $\a\in (1/2,1).$\smallskip

(b) Theorem 2 in \cite{NS} states that
$$M_\a\; =\; 1\; +\; \varepsilon \log(\varepsilon)\; +\; c_0\varepsilon \; +\; O(\varepsilon^2\log(\varepsilon))$$
as $\a\to 1,$ with $\varepsilon = 1-\a$ and $c_0 \sim -0.2228.$ This estimate is smaller than the one we get from (\ref{Best}), which is
$$1\; +\; \varepsilon \log(\varepsilon)\; +\; c_1\varepsilon \; +\; O(\varepsilon^2\log(\varepsilon))$$
with the same notation and $c_1= \log(4m_{1/2}) -1 \sim 0.4806.$ This shows that the median-mode inequality for $\Za$ holds also true as soon as $\a$ is close enough to 1. We believe that it is true for all $\a$'s.} 
\end{REMS}

Let us conclude this paragraph with partial mean-median-mode or median-mean inequalities for the Mittag-Leffler distribution. Set ${\tilde M}_\a, {\tilde m}_\a, {\tilde \mu}_\a$ for the respective mode, median and mean of $\M_\a$ and recall that $ {\tilde \mu}_\a =1/\Ga(1+\a).$ It is known that $\M_\a$ is always strictly unimodal and that ${\tilde M}_\a =0$ if and only if $\a\le 1/2$ - see Theorem (b) in \cite{TS3}. By Lemma 1.9 and Theorem 1.14 in \cite{DJ}, this shows that
\begin{equation}
\label{3MML}
{\tilde M}_\a \; < \; {\tilde m}_\a\; <\; {\tilde \mu}_\a
\end{equation}
for every $\a\in[0,1/2].$ The following proposition shows that (\ref{3MML}) is, however, not always true.

\begin{PROP} 
\label{4MML}
With the above notation, one has
$${\tilde m}_\a\; > \; {\tilde \mu}_\a$$
as soon as $\a \in [1-m_\S,1).$
\end{PROP}

\proof From the previous discussion, we have
$${\tilde m}_\a\; =\; m_\a^{-\a}\; \ge\; \a^{-\a}\; >\; \frac{1}{\Ga(1+\a)} \; =\; {\tilde \mu}_\a$$
for every $\a \in[1-m_\S,1)$ the strict inequality following from a direct analysis on the Gamma function.

\endproof

\begin{REMS}{\em (a) An upper bound for ${\tilde M}_\a$ when $\a > 1/2,$ which is not very explicit, can be obtained from Example 2 p. 307 in \cite{SY}.

\smallskip

(b) When $\a > 1/2,$ the number ${\tilde M}_\a > 0$ is also the unique mode of $X_1,$ where $\{X_t, \; t\ge 0\}$ is the spectrally negative $(1/\a)-$stable process defined before the statement of Theorem B - see Exercise 29.7 in \cite{S}. Differentiating the Laplace transform yields $\EE[X_1] = 0$ and it is conjectured that the median of $X_1$ lies inside $(0,{\tilde M}_\a),$ in other words that the mean-median-mode inequality holds for $X_1$ in reverse order. In general, the validity of the mean-median-mode inequality in one or another direction for {\em all} real $\a-$stable densities ($\a\in(1,2]$) is an open problem, which is believed to be challenging.

\smallskip

(c) It does not seem that the mean-median-mode inequality holds true in one or another direction for all $\M_\a$'s. Indeed,  if this were the case, then by (\ref{3MML}), Proposition \ref{4MML} and the obvious continuity of the three parameters with respect to $\a$ there would exist $\a\in (1/2, 1)$ such that 
$${\tilde M}_\a \; = \; {\tilde m}_\a\; =\; {\tilde \mu}_\a.$$
However the empirical fact, which can be observed e.g. on the graphics of \cite{Z1} p. 146, that the above $X_1$'s are completely skew to the left, in other words that the branch of their density on the left-hand side of the mode is everywhere above the branch on the right-hand side, would then lead to 
$$(\M_\a - {\tilde m}_\a)_- \;\elaw\; (\M_\a - {\tilde M}_\a)_-\; \prst\;(\M_\a - {\tilde m}_\a)_+ \;\elaw\; (\M_\a - {\tilde M}_\a)_+,$$
with the notation of \cite{DJ} p. 34. By the proof of Theorem 1.14 in \cite{DJ} this would entail ${\tilde \mu}_\a >  {\tilde m}_\a,$ a contradiction.}
 \end{REMS}

\subsection{Uniform estimates for the Mittag-Leffler function}

The classical Mittag-Leffler function is defined by
\begin{equation}
\label{Mittag}
\Ea (z) \; =\; \sum_{n\ge 0} \frac{z^n}{\Ga (1+\a n)}, \quad z \in \CC,
\end{equation}
and is a natural extension of the exponential function. Since its introduction in \cite{ML} it has been the study of many studies - see Chapter XVIII in \cite{E3} for classical properties and \cite{HMS} for a recent account with updated references. In this paragraph we are interested in the function 
$$x\mapsto \Ea(-x), \quad x\in\rl^+,$$ 
which is completely monotonic for any $0< \a < 1$ in view of the following Laplace representation (see Chapter XIII.8 (b) in \cite{F} and the references therein):
\begin{equation}
\label{Polar2}
\Ea(-x)\; =\; \EE[e^{-x\M_\a}], \quad x\ge 0.
\end{equation}
The asymptotic expansion at infinity for this function can be found e.g. in (6.11) of \cite{HMS}, and yields at the first order
\begin{equation}
\label{Mittaginf}
\Ea(-x)\; =\; \frac{1}{\Ga(1-\a) x}\, +\, o(x^{-1}), \quad x\to +\infty.
\end{equation}
Consider now the decreasing family of functions
$$f_c(x)\; =\; \frac{1}{1+cx}, \quad x\ge 0,$$
indexed by $c >0.$ This family converges uniformly on compact sets of $(0, +\infty)$ to $0$ when $c\to+\infty$ resp. to $1$ when $c\to 0.$ Combining this with the first-order expansions of $f_c$ and $\Ea$ at zero and infinity, one easily deduces the existence of $0 < c_1 < c_2 < +\infty$ such that 
$$f_{c_2}(x)\;\le \; \Ea(-x)\;\le \; f_{c_1}(x)$$
for all $x\in\rl^+.$ Defining $c_\a^- = \max\{c > 0\,\slash\, f_c(x)\ge  \Ea(-x),\, x\in\rl^+\}$
and $c_\a^+ = \min\{c > 0\,\slash\,\Ea(-x)\ge f_c(x),\, x\in\rl^+\}$ for every $\a\in (0,1),$ it is easy to see from (\ref{Mittag}), (\ref{Mittaginf}) and the first order expansions of $f_c$ at zero and infinity that necessarily one has
$$c_\a^-\; \le \; \Ga(1+\a)^{-1}\qquad\mbox{and}\qquad c_\a^+\; \ge \; \Ga(1-\a).$$
From numerical simulations, it is conjectured in Section 3.2 of \cite{M} - see especially (3.11) therein - that the above inequalities are actually equalities for every $\a\in (0,1).$ Notice in passing that the extension of these optimal bounds to $\a =1$ is simply 
\begin{equation}
\label{Fact}
0\;\le\; e^{-x}\;\le \;\frac{1}{1+x}
\end{equation}
for all $x\in\rl^+,$ a fact which is certainly true. A consequence of Theorem C is the validity of the above conjecture.
\begin{THE} 
\label{Maynard}
With the above notations, one has $c_\a^-= \Ga(1+\a)^{-1}$ and $c_\a^+ =\Ga(1-\a).$ In other words, for every $\a\in (0,1)$ the uniform estimate
\begin{equation}
\label{Main}
\frac{1}{1+\Ga(1-\a)x}\;\le \; \Ea(-x)\; \le\;\frac{1}{1+\Ga(1+\a)^{-1} x}
\end{equation}
holds over $\rl^+,$ with optimal constants.
\end{THE}

\proof To obtain the first inequality in (\ref{Main}), we write
\begin{eqnarray*}
\frac{1}{1+\Ga(1-\a) x}\; =\; \EE\lcr e^{-x\Ga(1-\a)\L}\rcr & = & x\int_0^\infty e^{-xt}\, \pb[\Ga(1-\a)\L\le t]\, dt\\
& \le & x\int_0^\infty e^{-xt}\, \pb[\M_\a \le t]\, dt\; =\; \EE\lcr e^{-x \M_\a}\rcr\; =\; \Ea(-x)
\end{eqnarray*}
for every $x > 0,$ where the second equality is an integration by parts, the inequality comes from $\M_\a\,\prost\,\Ga(1-\a) \L,$ and the last equality is (\ref{Polar2}). The optimality of $\Ga(1-\a)$ in (\ref{Main}) was already explained above and the validity of the inequality for $x = 0$ is trivial.

The second inequality in (\ref{Main}) is an analogous consequence of Theorem C. Integrating twice by parts, we get
\begin{eqnarray*}
\frac{1}{1+\Ga(1+\a)^{-1} x}\; =\; \EE\lcr e^{-x\Ga(1+\a)^{-1}\L}\rcr & = & x^2 \int_0^\infty e^{-xt}\lpa\int_0^t \pb[\Ga(1+\a)^{-1}\L\le s]\, ds\rpa dt\\
& \ge &  x^2\int_0^\infty e^{-xt}\lpa\int_0^t \pb[\M_\a \le s]\, ds\rpa dt\; =\;  \Ea(-x)
\end{eqnarray*}
for every $x > 0,$ where the inequality comes from $\Ga(1+\a)\M_\a\,\prcvx\,\L.$ The optimality of $\Ga(1+\a)$ in (\ref{Main}) was already discussed above and the validity of the inequality for $x = 0$ is obvious.

\endproof

\begin{REMS} {\em (a) When $\a\to 0+$ the above estimate entails that 
$$0\; \le\; \frac{1}{1+x} - \Ea(-x)\; \le \;\Ga(1-\a) - 1\; = \;\a\gamma\; +\; o(\a^2)$$ 
uniformly on $\rl^+,$ where $\gamma = - \Ga'(1)$ is Euler's constant. It does not seem possible to deduce this upper bound, which basically controls the speed of the convergence in law (\ref{Cresson}), from the mere definition of $\Ea$.\smallskip

(b) It follows from the above proof that
$$\frac{1}{x}\lpa \Ea(-x) - \frac{1}{1+\Ga(1-\a) x}\rpa\quad \mbox{and} \quad \frac{1}{x^2}\lpa \frac{1}{1+\Ga(1+\a)^{-1} x} - \Ea(-x) \rpa$$
are completely monotonic functions. \smallskip

(c) It follows from Theorem B, Theorem \ref{Mike} and the second part of the above proof that for all $x\ge 0$ 
\begin{equation}
\label{F2}
e^{-x} \; \le \; \Ea(-\Ga(1+\a)x)\; \le \; E_\bt(-\Ga(1+\bt)x)\; \le\; \frac{1}{1+x}
\end{equation} 
if $1/2\le \bt<\a <1$ or $\bt \le 1/2$ and $\a \ge 1/2\vee(2\bt\wedge (\bt+1)/2).$ This can be viewed as an improvement on (\ref{Fact}). From Remark \ref{F1} we believe that (\ref{F2}) remains true for all $0< \bt<\a <1.$ Observe that this would entail that the function $\a\mapsto \Ea(-x)$ is non-increasing on $[0,\alpha_0)$ for every $x\ge 0,$ where $1+\alpha_0 \sim 1.46163$ is the location of the minimum of the Gamma function. \smallskip

(d) It follows from Theorem \ref{Mike} and the first part of the above proof that for all $x\ge 0$
\begin{equation}
\label{F3}
 \frac{1}{1+x} \; \le \; E_\bt(-\Ga(1-\bt)^{-1}x)\; \le\; \Ea(-\Ga(1-\a)^{-1}x)\; \le \;1
\end{equation} 
if $1/2\le \bt<\a <1$ or $\bt \le 1/2$ and $\a \ge 1/2\vee(2\bt\wedge (\bt+1)/2).$ Again, we believe that (\ref{F3}) remains true for all $0< \bt<\a <1.$ Putting (\ref{F2}) and (\ref{F3}) together would draw the complete picture around Theorem \ref{Maynard}. 
}
\end{REMS}

\section{Factorizing one-sided stable branches}

Consider $X(\a,\rho)$ a strictly $\a-$stable random variable on the line with positivity parameter $\rho = \pb[X (\a,\rho) \ge 0],$ and normalized such that
$$\log[\EE[e^{\i \lbd X(\a,\rho)}]]\; =\; -(\i \lbd)^\a e^{-\i\pi\a\rho\, {\rm sgn}(\lbd)}, \quad \lbd\in\rl.$$
Recall that with this parametrisation which is (C) in the introduction of \cite{Z1} or (3) in \cite{Z2} p.355, one has $\rho\in[1-1/\a, 1/\a]$ if $\a\in (1,2]$ and $\rho\in[0,1]$ if $\a\in (0,1],$ with $X(1,1) = \Un$ and $X(1,0) = -\Un.$ We will implicitly exclude the two latter and only deterministic situations in the sequel. Comparing with Parametrisation (B) in \cite{Z1} or (1) in \cite{Z2} yields the following more familiar form for the characteristic exponent of $X(\a,\rho):$ 
$$\log[\EE[e^{\i \lbd X(\a,\rho)}]]\; =\;c\,\vert\lbd\vert^\a (1 - \i\theta\tan(\frac{\pi\a}{2})\,{\rm sgn}(\lbd))$$
with $\rho = 1/2 + (1/\pi\a) \tan^{-1}(\theta\tan(\pi\a/2))$ and $c = \cos(\pi\a(\rho -1/2)).$ Introduce the positive one-sided branch of $X(\a,\rho)$ which is the conditioned random variable
$$X^+(\a,\rho)\; =\; X(\a, \rho)\;\vert\; X(\a, \rho) \, \ge\, 0.$$
Observe that since $X(\a, \rho)\elaw -X(\a, 1-\rho),$ the density of $X (\a,\rho)$ can be recovered in pasting those of $X^+(\a,\rho)$ and $X^+(\a,1-\rho).$ This explains the natural interest in the random variable $X^+(\a,\rho),$ which is the matter of Chapter 3 in \cite{Z1}. 

In this section we would like to point out an extension of the factorizations stated in Theorem 1 to the random variable $X^+(\a,\rho),$ in the case when $\a,\rho$ are rational. It is a consequence of the following important proposition, which is known from  Formula (3.3.16) in \cite{Z1} but we provide a separate proof for the sake of clarity.

\begin{PROP}
\label{Folk}
With the above notation, one has
\begin{equation}
\label{Volker}
X^+(\a,\rho)\; \elaw\;\lpa\frac{\Z_{\a\rho}}{\Z_\rho}\rpa^\rho.
\end{equation}
\end{PROP}

\proof Suppose that $\rho\ge 1/2.$ Bochner's subordination for L\'evy processes - see Chapter 6 in \cite{S} for details, self-similarity and Exercise 29.7 in \cite{S} entail
$$X^+(\a,\rho)\; \elaw\;X^+(1/\rho,\rho)\times\Z_{\a\rho}^\rho\; \elaw\;\M_\rho\times\Z_{\a\rho}^\rho\; \elaw\;\lpa\frac{\Z_{\a\rho}}{\Z_\rho}\rpa^\rho.$$
On the other hand, setting $\rho' = 1/\rho-1$ and using Zolotarev's duality - see (2.3.3) in \cite{Z1} - we have
\begin{eqnarray*}
X^+(\a,1-\rho)\; \elaw\;X^+(1/\rho,1-\rho)\times\Z_{\a\rho}^\rho & \elaw & X^+(\rho,\rho')^{-\rho}\times\Z_{\a\rho}^\rho \\
&\elaw & X^+(1/\rho',\rho')^{-\rho}\times\Z_{\rho\rho'}^{-\rho\rho'} \times\Z_{\a\rho}^\rho \\
&\elaw & \lpa\frac{\Z_{\rho'}\times\Z_{\a\rho}^{\frac{1}{\rho'}}}{\Z_{\rho\rho'}}\rpa^{\rho\rho'} \;\elaw \; \lpa\frac{\Z_{\a(1-\rho)}}{\Z_{1-\rho}}\rpa^{1-\rho}.
\end{eqnarray*}
This completes the proof.
\endproof

\begin{REMS} {\em Let $\{X_t, \, t\ge 0\}$ be the spectrally negative $(1/\a)-$stable process introduced before the statement of Theorem B. The above proposition entails with self-explanatory notations that
$$X_1^- \;\elaw\; \M_{1-\a}\,\times\, \M_{\frac{1}{\a}-1}^{-\a}\; \rightarrow \; \frac{\L}{\L}\;\; \mbox{as $\a\to 1.$}$$
The limit in law has been observed in Formula (3.1.4) of \cite{Z1}.}
\end{REMS}

It is clear that putting (\ref{BetaGamma}) and (\ref{Volker}) resp. (\ref{Beta}) and (\ref{Volker}) together, we get a factorization of $X^+(\a,\rho)$ with $\a,\rho$ rational in terms of Beta and Gamma random variables resp. in terms of Beta random variables only. An interesting feature of these and all factorizations pertaining to (\ref{Volker}) is that they have a more down-to-earth character than the analytical representations in terms of Meijer's $G-$function - see the main result of \cite{Z2} p. 358 and also \cite{HP, H} for  essentially equivalent formulations. As an example, let us mention that (\ref{Volker}) can be used to provide a very quick proof of the unimodality of all stable densities \cite{TS2}.\smallskip

It is well-known that factorizations analogous to (\ref{Volker}) can also be obtained for the norm of isotropic $\a-$stable vectors. Consider on $\rl^d$ endowed with the standard Euclidean structure the random vector $\X_\a\, (0<\a\le 2)$ having characteristic function
$$\log[\EE[e^{\i <z, \X_\a>}]]\; =\;-\,\vert\vert z\vert\vert^\a.$$
When $\a=2$ one has $\X_2\elaw\cN(0, \sqrt{2}\,{\rm Id})$ and when $\a<2,$ Bochner's subordination entails
$$\vert\vert\X_\a\vert\vert\;\elaw \;\vert\vert\X_2\vert\vert\,\times\, \Z_{\frac{\a}{2}}\;\elaw \; 2\sqrt{\G_{d/2}}\,\times\, \Z_{\frac{\a}{2}}.$$
It would be interesting to see if this identity and the results of the present paper could not provide any further distributional properties of $\vert\vert\X_\a\vert\vert.$ As a first basic example, we can prove that the random variable $\vert\vert\X_\a\vert\vert^s$ is unimodal for all $s >0$ (with a positive mode if $d\ge 3$). \smallskip

We would like to conclude this paper with an open question on real stable densities which is related to Theorem B. With the above notation, this latter result states that the map $\rho \mapsto \Ga(1+\rho)X^+(1/\rho, \rho)$ is non-increasing for the convex order on $[1/2,1].$ Plotting the set of admissible parameters of $X^+(\a, \rho)$ with $\rho$ on the $x-$coordinate and $\a$ on the $y-$coordinate yields the sketch of a house with a roof of the Asian type. With this sketch in mind, Theorem B shows that moving along the top of the roof from left to right we get on $[1/2,1]$ a decreasing family for the convex order after suitable normalization. The following proposition shows that this phenomenon also occurs when moving along the lintel.

\begin{PROP}
\label{Lintel}
With the above notation, the map $\rho\mapsto \Ga(1+\rho)\Ga(1-\rho)X^+(1,\rho)$ is non-increasing for the convex order on $(0,1).$
\end{PROP}

\proof Since the involved random variables have no expectation, we only have to consider non-increasing convex functions in (\ref{cxx}). Setting $c_\rho =  \Ga(1+\rho)\Ga(1-\rho)$ for simplicity, this amounts to show that the map $\rho\mapsto c_\rho X^+(1,\rho)$ is non-decreasing for the increasing concave order, with the notation of Chapter 3 in \cite{SS}. This is an obvious consequence of the following claim
$$c_{\rho'}X^+(1,\rho')\;\prost\; c_{\rho} X^+(1,\rho)\qquad\mbox{if $\rho'\le\rho.$}$$
The latter claim is obtained in considering the explicit density of $c_\rho X^+(1,\rho),$ which is given by
$$\frac{1}{x^2 + 2c_\rho\cos(\pi\rho)x +c_\rho^2}$$
on $\rl^+,$ and using the fact that $\rho\mapsto c_\rho$ increases on $(0,1).$ 

\endproof

In the above proposition the convex order is not the most appropriate one because of the absence of expectations. But this result must be viewed as a limit of what should happen inside the roof: the conjecture is that for every $\a\in (1,2),$ there exists an Asian subroof connecting the parameters $(1,0), (\a, 1/2)$ and $(1,1)$ such that moving along the top of this subroof from left to right, we get a family of one-sided stable branches which is non-increasing for the convex order after suitable normalization. 

\bigskip

\noindent
{\bf Acknowledgements.} I thank Claude Lef\`evre for several discussions on stochastic orders at Universit\'e Libre de Bruxelles, and a referee for some useful comments. Ce travail a b\'en\'efici\'e d'une aide de l'Agence Nationale de la Recherche portant la r\'ef\'erence ANR-09-BLAN-0084-01.

\end{document}